\documentclass[11pt]{article}
\usepackage{graphicx,psfrag}
\usepackage{amssymb,amsfonts,amsthm}
\usepackage{amsmath,amsthm,amssymb}
\usepackage{amsfonts}
\usepackage{relsize}
\usepackage[T1]{fontenc}
\usepackage[utf8]{inputenc}
\usepackage{mathtools}   
\usepackage{graphicx}
\usepackage[OT2,OT1]{fontenc}
\newcommand\cyr{%
\renewcommand\rmdefault{wncyr}%
\renewcommand\sfdefault{wncyss}%
\renewcommand\encodingdefault{OT2}%
\normalfont
\selectfont}
\DeclareTextFontCommand{\textcyr}{\cyr}

\newcommand{\ab}{{\bar a}}

\newcommand{\bb}{{\bar b}}

\newcommand{\inp}{{\mathrm{inp}}}
\newcommand{\Sym}{{\mathrm{Sym}}}
\newcommand{\Add}{{\mathrm{Add}}}
\newcommand{\Sub}{{\mathrm{Sub}}}

\newcommand{\mor}{\orr{\mathcal M}}
\newcommand{\mol}{\oll{\mathcal M}}
\newcommand{\mm}{\mathcal M}
\newcommand{\vv}{\mathcal V}

\newcommand{\fin}{\mathrm{fin}}

\newcommand{\la}{\langle}
\newcommand{\ra}{\rangle}

\newtheorem{theorem}{Theorem}[section]
\newtheorem{lemma}[theorem]{Lemma}

\theoremstyle{definition}

\newtheorem{remark}[theorem]{Remark}

\newcommand {\N}{\mathbb{N}} 
\newcommand {\Z}{\mathbb{Z}}            

\newcommand {\me}{\medskip}
\newcommand{\oll}{\overleftarrow}
\newcommand{\orr}{\overrightarrow}

\newcommand {\iv}{^{-1}}

\begin{document}

\title{Minsky machines and algorithmic problems}

\author{Mark Sapir\thanks{The
research was supported in part by the NSF grants DMS 1418506 and
DMS	1318716, and a BSF grant.}}

\date{}
\maketitle

\me
%
%

\begin{abstract} This is a survey of using Minsky machines to study algorithmic problems in semigroups, groups and other algebraic systems.
\end{abstract}
\tableofcontents

\section{Introduction}

In 1966, Yuri Gurevich \cite{Gur} proved that the universal theory of finite semigroups is undecidable. One can interpret that result in several ways. For example, it means that given a number of semigroup relations $u_i=v_i$ and another relation $u=v$, we cannot algorithmically decide if the equality $u=v$ holds in every finite semigroup satisfying all the relations $u_i=v_i$. In that sense the universal theory of finite semigroups can be called the \emph{uniform word problem} of finite semigroups. Note that individually every finite semigroup has, of course, decidable word problem. Gurevich's result means that there is no uniform algorithm that works for all finite semigroups. That result turned out to be influential for two reasons. First, it opened the area of studying the uniform word problem in several classes of algebras, including semigroups and groups. Second, that was one of the first applications of Minsky (register) machines in proving undecidability of an algorithmic problem in algebra (the first result was the proof of undecidability of exponential diophantine equations from \cite{DPR}). The goal of this paper is to survey some applications of Minsky machines to various algorithmic problems in semigroups, groups and other types of algebras. Note that there is some intersection of this paper with the (250-page) survey paper \cite{KS} about algorithmic problems in varieties. But most  results surveyed here are not about varieties and are more recent than \cite{KS}.

\section{Turing machines and Minsky machines}

\subsection{Turing machines}
In this paper, we shall consider several types of machines. A machine $M$ in general has an alphabet and a set of
words in that alphabet called  configurations (by ``words" we may also mean ``numbers written in binary or unary" or ``tuples of numbers"). It also has a finite set of commands. Each command is a
partial injective transformation of the set of configurations. A machine  is called {\em deterministic} if the domains of its commands are disjoint.
A machine usually has a distinguished {\em stop} configuration, and a set $I=I(M)$ of {\em input} configurations.

A {\em computation} of $M$ is a finite or infinite sequence of configurations and commands from $P$:

$$w_1 \stackrel{\theta_1}{\longrightarrow} w_2 \stackrel{\theta_2}{\longrightarrow} \ldots
\stackrel{\theta_l}{\longrightarrow} w_{l+1},\ldots $$
such that $\theta _i (w_i)=w_{i+1}$ for every $i=1,\ldots ,l,\ldots .$

If the computation is finite and $w_{l+1}$ is the last configuration, then  $l$ is called the {\em length} of the computation. A
configuration is called {\em accepted} by $M$ if there exists a computation connecting that configuration with
the stop configuration. The {\em time function} $T_M(n)$ of $M$ is the minimal function such that every accepted word of length $\le n$
has an accepting computation of length $\le T_M(n)$.

The machine $\Sym(M)$ is made from $M$ by adding the inverses of all commands of
$M$. Two configurations $w$, $w'$ are called {\em equivalent}, written $w\equiv_M w'$, if there exists a
computation of $\Sym(M)$ connecting these configurations. Clearly, $\equiv_M$ is an equivalence relation.

The following general lemma is an easy exercise but it is very useful.

\begin{lemma}\label{l:genmachine}\label{l:sym} Suppose that $M$ is deterministic. Then two configurations $w,w'$
of $M$ are equivalent if and only if there exist two computations of $M$ connecting $w,w'$ with the same
configuration $w''$ of $M$.
\end{lemma}
%

We say that a set $X$ of natural numbers is {\em enumerated} by a machine $M$ if there exists a recursive
encoding $\mu$ of natural numbers by input configurations of $M$ such that a number $u$ belongs to $X$ if and
only if $\mu(u)$ is accepted by $M$. The set $X$ is {\em recognized} by $M$ if $M$ enumerates $X$ and for every
input configuration every computation starting with that configuration, eventually halts  (arrives to a
configuration to which no command of $M$ is applicable).

We say that machine $M'$ {\em polynomially reduces} to a machine $M$ if there exists a {{deterministic}} polynomial time algorithm
$A$ checking equivalence of configurations of $M'$ which uses an oracle checking equivalence of configurations of
$M$.

We say that $M$ and $M'$ are {\em polynomially equivalent} if there are polynomial reductions of $M$ to $M'$ and
vice versa.

For example, a Turing machine $M$ with $K$ tapes
consists of hardware (the tape alphabet $A=\sqcup_{i=1}^k A_i$, and the state alphabet $Q=\sqcup_{i=1}^K
Q_i$\footnote{$\sqcup$ denotes disjoint union}) and program $P$ (the list of commands, defined below). A {\em
configuration} of a Turing machine $M$ is a word $$\alpha_1 u_1q_1v_1\omega_1\hbox{  } \alpha_2u_2q_2v_2\omega_2\hbox{  }
\ldots\hbox{  } \alpha_Ku_Kq_Kv_K\omega_K$$ (we included spaces to make the word more readable) where $u_i, v_i$ are
words in $A_i$, $q_i\in Q_i$ and $\alpha_i,\omega_i$ are special symbols (not from $A\cup Q$).  {{A tape of the machine is a part of the configuration, it is a subword
from $\alpha _i$ to $\omega _i$.}}

A command simultaneously replaces subwords $a_iq_ib_i$ by words $a_i'q_i'b_i'$ where $a_i,a_i',$ are either
letters from $A_i\cup\{\alpha_i\}$ or empty, $b_i,b_i'$ are either letters from $A_i\cup \{\omega_i\}$ or empty.
A command cannot insert or erase $\alpha_i$ or $\omega_i$, so if, say, $a_i=\alpha_i$, then $a_i'=\alpha_i$. Note
that with every command $\theta$ one can consider the {\em inverse} command $\theta\iv$ which undoes what
$\theta$ does.

For the Turing machine we choose {\em stop states} $q_i^0$ in each $Q_i$, then
 a configuration $w$ is accepted if there exists a computation starting with $w$ and ending with a
configuration where all state symbols are $q_i^0$ and all tapes are empty (which is the stop configuration for the Turing machine). Also we choose {\em start states}
$q_i^1$ in each $Q_i$. Then an  input configuration corresponding to a word $u$ over $A_1$ is a
configuration $\inp(u)$ of the  form $$\alpha_1uq_1^1\omega_1\hbox{  } \alpha_2q_2^1\omega_2\hbox{  } \ldots\hbox{  } \alpha_Kq_K^1\omega_K.$$
We say that a word $u$ over $A_1$ is accepted by $M$ if the configuration $\inp(u)$ is accepted. The set of all
words accepted by $M$ is called the {\em language accepted by} $M$.

\subsection{Minsky machines}

The hardware of a $K$-glass Minsky machine, $K\ge 2$, consists of $K$ glasses containing coins. We assume that
these glasses are of infinite height. The machine can add a coin to a glass, and remove a coin from a glass
(provided the glass is not empty). The number of coins in the glass $\# k$ is denoted by $\epsilon_k$.
%

In the \emph{program} of every Minsky machine, the commands are numbered, command $\# 1$ is the \emph{start command}, command $\# 0$ is the \emph{stop command}. A \emph{configuration} of a 2-glass Minsky machine is a triple of numbers $(i;m,n)$ where $i$ is the number of command being executed, $m$ is the number of coins in the first glass, $n$ is the number of coins in the second glass. The start configurations have the form $(1;m,n)$ and the stop configurations have the form $(0;m,n)$.

A {\em command} $\# i$, $i\ge 1$, has one of the following forms:

\begin{itemize}
\item Put a coin in the glass $\# k$ and go to command $\# j$. We shall encode this
    command by
$$\Add(k); j;$$
\item If the glass $\# k$ is not empty then take a coin from it and go to command $\# j$. This command is encoded by $$\epsilon_{k}>0 \to \Sub(k); j;$$
\item  If the glass $\# k$ is empty, then  go to instruction $\# j$. This command is encoded by
$$\epsilon_{k}=0\to j; $$
\end{itemize}

Note that here $j$ may be equal to $0$, but there is no instruction associated with command $\# 0$.


\begin{remark} This defines deterministic Minsky machines. We will also need non-deterministic Minsky machines.
Those
will have two or more commands with the same number.
\end{remark}

The proof of the following theorem can be found in Malcev \cite{malcev} or extracted from the proof there.

\begin{theorem}\label{t:MM}
Let $X$ be a  recursively enumerable set of natural numbers. Then the following holds:
\begin{itemize}
\item[(a)] there exists a 2-glass deterministic Minsky machine $M$ which recognizes
$L$ in the following sense: $M$   begins its work in configuration $(1;2^m,0)$ and halts in configuration
$(0;0,0)$ if  and only if $m\in X$,  and it  works forever  if $m\not\in X$.

\item[(b)] Every computation of $M$ starting with a configuration $c$
    empties each glass after at most $O(|c|)$ steps.


\end{itemize}
\end{theorem}

\section{The three main semigroups simulating Minsky machines}

There are three basic ways to interpret 2-glass Minsky machines in semigroups. They correspond to the three ways to put two glasses and the machine head (the one that counts the commands and puts coins in the glasses) on the line: the head can be between two glasses, to the left of the glasses and to the right of the glasses. If we imagine the head to have a short hand used to put the coins in glasses, then in the last two cases we should be able to permute the two glasses: if the counter wants to put a coin in, say, glass $\# 2$, and glass $\# 1$ is between the head and glass $\#2$, then first the two glasses are permuted.

The three semigroups corresponding to a 2-glass Minsky machine $M$ are $S_1(M)$, $\orr S_2(M)$ and $\oll S_2(M)$. The last two semigroups are anti-isomorphic, so we only define $S_1(M)$ an $\orr S_2(M)$. Each of the three semigroups has zero $0$.

Let $M$ be a Minsky machine with $2$
glasses and commands $\#\# 1,2,...,N, 0$. Then both semigroups $S_1(M)$ and $\orr S_2(M)$ is generated by the elements $q_0,q_1,\ldots,q_N$ and
$\{a_i,A_i, i=1,2\}$. Here $a_i$ play the role of coins in glass $\# i$, $q_i$ play the role of numbers of commands (i.e., the states of the head), and $A_i$ play the role of the bottoms of glasses (these are needed in order to be able to check if a glass is empty). The set of defining relations of $S_1(M)$ and $\orr S_2(M)$ contains the following relations corresponding to the commands of $M$.

$$
\begin{array}{|l|l|l|}
\hline
\hbox{Command \# i of } M, i\ge 1& \hbox{Relation of } S_1(M) & \hbox{Relation of } \orr S_2(M) \\
\hline
\Add(1);j & q_i=a_1q_j & q_i=q_ja_1\\
\hline
\Add(2);j & q_i=q_ja_2 & q_i=q_ja_2\\
\hline
\epsilon_{1}>0\to \Sub(1);j & a_1q_i=q_j & q_ia_1=q_j\\
\hline
\epsilon_{2}>0\to \Sub(2);j &q_ia_2=q_j & q_ia_2=q_j\\
\hline
\epsilon_1=0\to j & A_1q_i=A_1q_j & q_iA_1=q_jA_1\\
\hline
\epsilon_2=0\to j & q_iA_2=q_jA_2 & q_iA_2=q_iA_2\\
\hline
\end{array}
$$
These will be called the \emph{Minsky relations}. The semigroups $S_1(M), \orr S_2(M)$ also have the following \emph{auxiliary relation}.

{\bf The auxiliary relations of $S_1(M)$:}

\begin{itemize}
\item All 2-letter words in the generators of $S_1(M)$ that are not subwords of the words $A_1a_1^mq_ia_2^nA_2$, $m,n\in \N$, $i=1,...,N$, are equal to 0;

\item $q_0=0$.
\end{itemize}

{\bf The auxiliary relations of $\orr S_2(M)$:}

\begin{itemize}
\item (Glass permuting relations)  each letter in $\{a_1,A_1\}$ commutes with each letter in $\{a_2,A_2\}$;

\item All 2-letter words in generators of $\orr S_2(M)$ that are not subwords of the words that are equal to $q_ia_1^mA_1a_2^nA_2$ modulo the glass permuting relations are equal to 0;

\item $q_0=0$.
\end{itemize}

A configuration $(i;m,n)$ of the Minsky machine $M$ corresponds to the element $w_1(i;m,n)=A_1a_1^mq_ia_2^nA_2$ in $S_1(M)$ and to the element $w_2(i;m,n)=q_ia_1^mA_1a_2^nA_2$ in $\orr S_2(M)$.

\begin{remark} \label{r:1} The auxiliary relations ensure that every word in the generators of $S_1(M)$ (resp. $S_2(M)$) that is not equal to 0 is a subword of a word of the form $w_1(i;m,n)$ (resp. a subword of a word that is equal to one of the words  $w_2(i;m,n)$ modulo the glass permuting relations).
\end{remark}

\begin{remark} Thus the semigroups $S_1(M)$, $\orr S_2(M)$, $\oll S_2(M)$ basically consist of the subwords of the words corresponding to the configurations of $M$. This is a crucial property of Minsky machines which makes them much better suited for semigroup simulation than the general Turing machines.
\end{remark}

\section{Varieties of semigroups and the word problem}

The proof of the following statement is straightforward, it is discussed in \cite{KS}.

\begin{theorem}[See \cite{KS}] Let $M$ be a Minsky machine. For every two configurations $(i;m,n)$ and $(i';m',n')$ the words $w_s(i;m,n)$ and $w_s(i',m',n')$, $s=1,2,3$  are equal in $S_1(M)$ (resp. $\orr S_2(M)$ or $\oll S_2(M)$) if and only if the configurations are equivalent. In particular, if $M$ has undecidable halting problem, then the word problem in each of the three semigroups associated with $M$ is undecidable. Moreover, the equality to 0 is undecidable in these semigroups.
\end{theorem}

Of course constructing a finitely presented semigroup with undecidable word problem (first done by Markov and Post, see \cite{malcev}) is easy enough using the ordinary Turing machines. The advantage of Minsky machines is that the semigroups $S_1(M), \orr S_2(M), \oll S_2(M)$ are in some sense ``small".

For example while examples corresponding to the ordinary Turing machines usually contain non-commutative free subsemigroups, and hence do not satisfy non-trivial identities (laws),  it is easy to see that each of the three semigroups $S_1(M), \orr S_2(M), \oll S_2(M)$ satisfies a non-trivial identity. For example, each of them satisfies $x^2y^2=y^2x^2$. This follows immediately from Remark \ref{r:1} and the auxiliary defining relations of these semigroups. In fact, one can describe all identities satisfied by these semigroups \cite{SapirAlgProbVar}.

In particular, the following theorem holds. Let $e_{ij}$ denote the $2\times 2$-matrix unit $(i,j)$-entry 1 and all other entries $0$. Let $\orr P$ denote the three element semigroup $\{e_{11}, e_{12}, 0\}$,  $\oll P$ denote the three-element semigroup $\{e_{11},e_{21},0\}$ and $T$ denote the four element semigroup $\{e_{11}, e_{12}, e_{22},0\}$. For every semigroup $S$ let $S^1$ be the semigroup $S$ with identity element formally adjoint. For example $P^1$ is the four element semigroup $\{1,e_{11},e_{12},0\}$. Let $\N$ be the additive semigroup of natural numbers. It is well-known and trivial that an identity $u=v$ is true in $\N$ if and only if it is {\em balanced}, that is if every letter occurs the same number of times in $u$ and in $v$.

\begin{theorem}[Sapir, \cite{SapirAlgProbVar}] For every Minsky machine $M$ the variety generated by $S_1(M)$, i.e., the smallest class of semigroups containing $S_1(M)$ and given by identities coincides with the variety ${\mathcal M}_1$ generated by the direct product $T\times \N$. The variety generated by $\orr S_2(M)$ coincides with the variety $\orr{\mathcal M}_2$ generated by $\oll P\times \orr P^1\times \N$, and the variety generated by $\oll S_2(M)$ coincides with the variety $\oll{\mathcal M}_2$ generated by $\oll P^1\times \orr P\times \N$ (thus these varieties do not depend on $M$).
\end{theorem}

Moreover Minsky machines and the easy construction above proved to be the universal tool in dealing with the word problem in semigroups satisfying identities. In particular, one can completely describe non-periodic varieties containing finitely presented semigroups with undecidable word problem.

We say that a finitely generated semigroup $S$ is finitely presented inside a variety $\mathcal V$ if it is defined by the identities of $\mathcal V$ plus a finite number of relations. We shall need the following sequence of {\emph Zimin words}: $$Z_1=x_1,...,Z_{n+1}=Z_nx_{n+1}Z_n.$$ This sequence of words plays an important role in combinatorial algebra (see \cite{KS}, \cite{SapirBook}). We say that a word $W$ is {\emph not an isoterm for an identity} $u=v$ if for some substitution $\phi$ of words for letters of $u,v$ we have that $\phi(u)\ne \phi(v)$ but $W$ contains either $\phi(u)$ or $\phi(v)$ as a subword. For example, the word $ababbab$ is not an isoterm for the identity $x^2=x^3$ because the word contains a subword $abab$ which is equal to $\phi(x^2)$ under the substitution $x\to ab$, and $\phi(x^2)\ne \phi(x^3)$. But it is an isoterm for the identity $x^3=x^4$ because it does not contain any subword of the  form $u^3$ or $u^4$

\begin{theorem}[Sapir, \cite{SapirAlgProbVar}] \label{t:9} Let $\mathcal V$ be a variety of semigroups defined by identities $u_i=v_i$ in at most $n$ variables and non-periodic (i.e., containing the semigroup $\N$, or, equivalently, every identity $u_i=v_i$ is balanced).
Then the following conditions are equivalent.

\begin{enumerate}
\item[(1)] Every semigroup that is finitely presented inside $\mathcal V$ has decidable word problem.
\item[(2)] Every semigroup that is finitely presented inside $\mathcal V$ has decidable elementary theory.
\item[(3)] Every semigroup that is finitely presented inside $\mathcal V$ is faithfully representable by matrices over a field.
\item[(4)] The variety does not contain varieties $\mm_1, \mor_2, \mol_2$ and the word $Z_{n+1}$ is an isoterm for every identity $u_i=v_i$.
\end{enumerate}

\end{theorem}

As mentioned in \cite{SapirAlgProbVar}, Property (4) of Theorem \ref{t:9} is algorithmically verifiable given a finite number of identities $u_i=v_i$.

\section{Gurevich's theorem. The uniform word problem for finite semigroups}

Let $L$ be a finite conjunction of equalities $u=v$, where $u,v$ are words in some alphabet $X$. Let $U, V$ be two words in $X$. Then the universal formula $L\to U=V$ is called a quasi-identity.  We say that the {\em uniform word problem} is solvable in a class of semigroups $\mathcal V$ if there exists algorithm that, given a quasi-identity $\theta$ decides whether $\theta$ holds in $\mathcal V$. Clearly, the uniform word problem is solvable if $V$ consists of finitely many finite semigroups. For every variety $\mathcal V$ of semigroups $\mathcal V_{\fin}$ denotes the set of finite semigroups from $\mathcal V$. Yu. Gurevich proved \cite{Gur} that if $\mathcal V$ is the variety of all semigroups, then the uniform word problem is not decidable in $\mathcal V_{\fin}$. Generalizing that result, we completely described  in \cite{SapirAlgProbVar,SapirWeak} all finitely based varieties $\vv$ such that the uniform word problem is decidable in $\vv_{\fin}$.

\begin{theorem}[Sapir, \cite{SapirWeak}] \label{t:8} For every finite set of identities $\Sigma$ in $n$ variables defining a variety $\vv$ the following conditions are equivalent:

\begin{itemize}\item[(1)] The uniform word problem is decidable in $\vv_{\fin}$.
\item[(2)] The word $Z_{n+1}$ is not an isoterm for $\Sigma$ and either $\vv$ is periodic (i.e., does not contain $\N$, or, equivalently $\Sigma$ contains a non-balanced identity), or does not contain any of the semigroups $T$, $\orr P^1\times \oll P$ and $\orr P\times \oll P^1$.
\end{itemize}
\end{theorem}

The proof of Theorem \ref{t:8} proceeds as follows. $(1)\to (2)$. Suppose that the uniform word problem is decidable in $\vv_{\fin}$.

First suppose that $\vv$ contains one of the varieties $\mm_1, \orr \mm_2, \oll \mm_2$.  In order to get a contradiction, consider the following modification of semigroups $S_1(M), \orr S_2(M), \oll S_2(M)$. We add three letters $c, c', e, C$ to the generating set of each semigroup. The construction below uses ideas from \cite{Gur} and is somewhat easier than a construction from \cite{SapirAlgProbVar}. The Minsky relations of the semigroups $S_1'(M), \orr S_2'(M), \oll S_2'(M)$ are defined as follows.

$$
\begin{array}{|l|l|l|}
\hline
\hbox{Command \# i of } M, i\ge 1& \hbox{Relation of } S'_1(M) & \hbox{Relation of } \orr S'_2(M) \\
\hline
\Add(1);j & q_i=a_1q_j & q_i=q_ja_1\\
\hline
\Add(2);j & q_i=q_ja_2 & q_i=q_ja_2\\
\hline
\epsilon_1=1 \to \Sub(1);j &A_1a_1w_1=cA_1q_1 & q_1a_1A_1=q_1A_1c \\
\hline
\epsilon_{1}\ge 2\to \Sub(1);j & a_1q_i=q_j & q_ia_1=q_j\\
\hline
\epsilon_{2}>0\to \Sub(2);j &q_ia_2=q_j & q_ia_2=q_j\\
\hline
\epsilon_1=0\to j & A_1q_i=A_1q_j & q_iA_1=q_jA_1\\
\hline
\epsilon_2=0\to j & q_iA_2=q_jA_2 & q_iA_2=q_iA_2\\
\hline
\end{array}
$$

One can also view this modification as a modification of the machine $M$: we add a new \emph{counter} glass which ``counts" how many times during a computation the first glass becomes empty: every time we remove the only remaining coin from the first glass, we add a coin in the new glass.

{\bf The auxiliary relations of $S'_1(M)$} are
\begin{itemize}

\item (Counter relations) $cc'=c'c=e, ec=c=ce, ec'=c'e=c'$, that is the subsemigroup generated by $c,c'$ is a subgroup isomorphic to the group $\Z$ with identity element $e$.

\item Every 2-letter word in the generators of $S'_1(M)$ which is not a subword of a word that is equal to a word of the form $Cc^kA_1a_1^mq_ia_2^nA_2$ modulo the counter relations is equal to 0.

\item $CeA_1q_i=0$ for every $i$.

\end{itemize}

{\bf The auxiliary relations of $\orr S_2'(M)$} are

\begin{itemize}
\item (Counter relation) $cc'=c'c=e, ec=c=ce, ec'=c'e=c'$, that is the subsemigroup generated by $c,c'$ is a subgroup isomorphic to $\Z$.

\item (Glass permuting relations) Every letter from $\{a_1, A_1\}$ commutes with every letter from $\{a_2,A_2\}$, every letter from $\{a_1,a_2,A_1,A_2\}$ commutes with every letter from $\{c, c', e,C\}$

\item Every 2-letter word in the generators of $S'_1(M)$ which is not a subword of a word that is equal to a word of the form $q_i a_1^mA_1a_2^nA_2c^kC$ modulo the Counter and Glass permuting relations is equal to 0.
\end{itemize}

Now it is  proved in \cite{SapirAlgProbVar} that $S_1=S_1'(M)$ belongs to the variety $\mm_1$, and $S_2=\orr S_2'(M)$ belongs to the variety $\mm_2=\orr \mm_2$ generated by $\oll P\times \orr P^1\times \N$. Let $L_i$, $i=1,2$, be the conjunction of the defining relations of $S_i$.
For every input configuration $(1;m,n)$ of $M$ let $W_1(1;m,n)$ be the word $CA_1a_1^mq_1a_2^nA_2$ and $W_2(1;m,n)=q_1a_1^mA_1a_2^nA_2C$. Consider the quasi-identity $L_i\to W_i=0$. Suppose that a configuration $(1;2^m,0)$ is accepted by $M$. Then it is proved in \cite{SapirWeak} that there are only finitely many elements in $S_i$ that divide the element $W_i(1;2^m,0)$ (recall that we say that an element $b$ \emph{divides} element $a$ if for some $x_1,x_2$ we have $a=x_1bx_2$). Consider the set of all elements of $S_i$ that do not divide $W_i=W_i(1;2^m,0)$. This set is an ideal $J_i$ of $S_i$, not containing $W_i$. Consider the Rees quotient $F_i=S_i/J$. It is a finite semigroup and it belongs to $\mm_i$ as a quotient of $S_i$. We can assume that $F_i$ is generated by the same generating set as $S_i$. Therefore $F_i$ satisfies $L_i$. But $W_i \ne 0$ in $F_i$ because $W_i\not\in J$. Thus a finite semigroup in $\mm_i$ does not satisfy the quasi-identity $L_i\to W_i(1;2^m,0)=0$.

On the other hand, suppose that $M$ works indefinitely long starting with the configuration $(1;2^m,0)$. Then by Part (b) of Theorem \ref{t:MM} the word $W_i(1;2^m,0)$ is equal in $S_i$ to words of the form $Cc^sA_1q_ta_2^nA_2$ for arbitrary $s$. Let $R_i$ be a periodic semigroup satisfying the formula $L_i$ for some interpretation of its variables in $R_i$. In other words, let $R_i$ be a quotient of $S_i$ generated by the same generating set as $S_i$. Then $R_i$ satisfies the relations of the form $W_i(1;2^m,0)=Cc^sA_1q_ta_2^nA_2$ for every $s$. Since $R_i$ is periodic, $c^s=c^{2s}$ for some $s$. Since the subsemigroup generated by $c, c'$ in $S_i$ is a group (isomorphic to $\Z$) with identity element $e$, the subsemigroup generated by $c$ in $R_i$ is a finite cyclic group, and so $c^s=e$ for some $s$. But $CeA_1q_t=0$ in $S_1$, and hence in $R_i$. Therefore we can conclude that the quasi-identity $L_i\to W_i(1;2^m,0)=0$ holds in every periodic semigroup.

Since there exists a Minsky machine satisfying the conditions of Theorem \ref{t:MM} for which the language of accepted input configurations is not recursive, there is no algorithm separating the quasi-identities that hold in all finite semigroups of $\vv$ from the set of quasi-identities that hold in all periodic semigroups. Thus the uniform word problem in $\vv$ is undecidable, a contradiction.

Now suppose that $Z_{n+1}$ is an isoterm for identities from $\Sigma$. Then it is proved in \cite{SapirBurns} (see also \cite{SapirBook}) that there exists a finite alphabet $x_1,...,x_k$ and a substitution $\phi\colon x_i\mapsto \phi(x_i)$ where $\phi(x_i)$ is a word in $\{x_1,...,x_k\}$, such that $\phi^s(x_1)$ is an isoterm for every identity from $\Sigma$ for every $s$. In that case, we construct in \cite{SapirWeak} another semigroup $S(M)$ simulating arbitrary 2-glass Minsky machine. In that semigroup the number $s$ of coins in glass $\# j$ is simulated not by the power $a_j^s$ (as in the constructions above) but by the word $\phi^s(x_1)$. As a result the words corresponding to configurations of $M$ are isoterms of the identities of $\Sigma$, and the situation is  similar to the situation with the variety of all semigroups (since the identities of the variety cannot ``mix up" the words corresponding to configurations of the Minsky machine). Of course replacing $a_i^s$ by $\phi^s(x_1)$ costs us something. For example, to simulate one command of $M$, we need several relations of $S(M)$, but it can be done, see \cite{SapirWeak}. This proves that in the case when $Z_{n+1}$ is not an isoterm for $\Sigma$, the uniform word problem in $\vv$ is also undecidable, which concludes the proof of implication $(1)\to (2)$ in Theorem \ref{t:8}.

$(2)\to (1).$ Suppose that the conditions of (2) hold. Then either $\vv$ is periodic or it contains $\N$. If it is periodic, then we proved in \cite{SapirRestrBurns} that the \emph{restricted Burnside} property holds in $\vv$. This means that for every natural number $k\ge 1$ there are only finitely many (effectively computable) finite semigroups in $\vv$ with at most $k$ generators. To prove that, we were using the celebrated positive solution of the restricted Burnside problem for groups  by Zelmanov \cite{Ze90,Ze91}, and our results on Burnside problems in semigroup varieties \cite{SapirBurns}.  This property easily implies solvability of the uniform word problem in $\vv_{\fin}$.

If $\vv$ is not periodic and  contains none of the three semigroups listed in Part (2) of the theorem, then by Theorem \ref{t:9} every semigroup that is finitely presented in $\vv$ is faithfully represented by matrices over a field, hence residually finite \cite{Malcev}. This also implies solvability of the uniform word problem in $\vv_{\fin}$. Indeed, in order to check if a quasi-identity  $L\to W=W'$ holds in $\vv_{\fin}$, consider the semigroup $E$ defined by the relations from $L$ in $\vv$. If $W\ne W'$ in that semigroup, then $W\ne W'$ in some finite quotient $E'$ of $E$. Thus $E'$ does not satisfy $L\to W=W'$. Hence $L\to W=W'$ does not hold in $\vv_{\fin}$. On the other hand, if $W=W'$ holds in $E$, then this equality holds in every homomorphic image of $E$, hence $L\to W=W'$ holds in $\vv_{\fin}$. Since the word problem is decidable in $E$ by Theorem \ref{t:9}, we can decide whether or not $W=W'$ in $E$, and hence whether or not $L\to W=W'$ holds in $\vv_{\fin}$.

\section{The uniform word problem for finite groups}

\subsection{Slobodskoi's theorem}

The first simulation of Minsky machines in groups was done by Slobodskoi \cite{sl}. This is not as easy as in the case of semigroups. The main problem is that we cannot simulate a command of a Minsky machine, say, $$i: \Add(1); j$$ by a substitution $q_i=q_ja_1$ because then the group would collapse. For example, if $i=j$, then we would have $q_i=q_ia_1$ and $a_1=1$. Thus the idea is to abandon the product operation in the group and use some other {\em derived} operation, say, commutator. Thus if $*$ is the new binary operation, we can simulate the command of a Minsky machine by $q_i=q_j*a_i$ and since $*$ does not necessarily satisfies the cancellation property, we avoid at least immediate collapse of the group. But then several new problems occur. For example $*$ will not be associative, and so we would have to interpret a configuration $(i;m,n)$ of a Minsky machine by a nested word like $(...(q_i*a_1)...)*a_1)*A_1)...$. In this case, it may be difficult to simulate permutation of glasses. These difficulties and ways to resolve them are described in details in \cite{KS}. Slobodskoi proved

\begin{theorem}[Slobodskoi, \cite{sl}] \label{t:sl} The uniform word problem us undecidable in the class of finite groups.
\end{theorem}

Several generalizations of Theorem \ref{t:sl} were then proved by Kharlampovich (see \cite{Khnil,Khdiss}), each time Minsky machines were used. Finally in \cite{KMS}, we proved the following result, again using Minsky machines.

\begin{theorem}\label{t:sl} Let ${\mathcal G}_1$ be the set of finite groups $G$ which have a normal series $N_1\unlhd N_2\unlhd G$ such that $N_1$ and $N_2$ are Abelian groups of the same prime exponent $p$, and $G/N_2$ is Abelian. Let ${\mathcal G}_2$ be the set of finite groups $G$ with normal series $N_1\unlhd N_2\unlhd G$ such that $N_1$ is contained in the center of $G$, $N_2/N_1$ is nilpotent of class at most 5, and $G/N_2$ is Abelian. Let ${\mathcal G}={\mathcal G}_1\cap {\mathcal G}_2$. Then the uniform word problem is undecidable in any set of finite groups containing ${\mathcal G}$.
\end{theorem}

\subsection{Some applications of Slobodskoi's theorem}

\subsubsection{Embeddings of finite semigroups into simple semigroups}

A semigroup with 0 is called \emph{0-simple} if it does not have any ideals except $0$ and itself (i.e., every two non-zero elements divide each other). These semigroups play very important role in the theory of semigroups being the building blocks from which all other semigroups are constructed. Finite 0-simple semigroups have very explicit structure (proved by Sushkevich and Rees independently \cite{CP}). For every such semigroup $S$ there exists a (finite) group $G$ and a $m\times n$-matrix $P$ where every entry is an element of $G$ or 0. The elements of $S$ are 0 and all triples $(i,g,j)$ where $1\le i\le m, 1\le j\le n, g\in G$. The product $(i,g,j)(i',g',j')$ is defined as 0 if $P(i',j)=0$ and $(i,gP(i',j)g',j')$ if $P(i',j)\ne 0$. Since finite 0-simple semigroups are so easy, it would be natural to guess that the set of their subsemigroups is also easy. That was proved not to be the case by Kublanovsky (first published in \cite{HKMST}). For every finite partial group $G$, that is a finite set with a partial operation $\cdot$ he constructed a finite effectively computable set of 4-nilpotent finite semigroup $N_i(G)$. This set satisfies the property that $G$ embeds into a group if and only if one of the $N_i(G)$ embeds into a finite 0-simple semigroup. The semigroups $N_i(G)$ are constructed as follows. We can assume $G$ contains the identity element, every element of $G$ has an inverse, and for every $a,b,c\in G$ we have $(a\cdot b)\cdot a\cdot (b\cdot c)$ provided each of the products involved in that equality is defined. Indeed, otherwise $G$ cannot be embeddable into a group. Consider all (finitely many) partial groups $G_i$ such that

\begin{itemize}

\item $G\le G_i$;

\item $G_i=G\cdot G \cdot G$.
\end{itemize}

Then $N_i(G)=N(G_i)$ is the semigroup defined on the set

$$(\{1\} \times  G \times \{2\}) \cup (\{2\} \times G\times  \{3\}) \cup (\{3\} \times G \times \{4\})
\cup (\{1\} \times (G\cdot G) \times \{3\}) \cup (\{2\} \times (G\cdot G) \times \{4\})
\cup G_i\cup \{0\}$$
with the operation $(i,u,j)(j,v,k)=(i,u\cdot v,k)$ if the right hand side is in $N_i(G)\setminus \{0\}$ or $0$ otherwise (that construction uses the idea of \emph{split systems} from my undergraduate diploma thesis \cite{SapirDiploma}).

It is well known \cite{Evans} that the problem whether a finite partial group embeds into a finite group is decidable if and only if the uniform word problem for finite groups is decidable. Thus Slobodskoi's Theorem \ref{t:sl} implies

\begin{theorem}[Kublanovsky \cite{HKMST}] \label{t:ku} The set of (4-nipotent) subsemigroups of finite 0-simple semigroups is not recursive.
\end{theorem}

\subsubsection{Equations over finite semigroups, the Rhodes' problem}

Let $S$ be a finite semigroup, $X=\{x_1,...,x_n\}$ be a set of variables. Let $u_i,v_i$ be words in $S\cup X$. Then the set of equalities $u_i=v_i$ is called a \emph{system of equations over} $S$. We say that a system of equations over $S$ is \emph{solvable} if there exists a finite semigroup $T\ge S$ and a map $X\to T$ which makes all equalities $u_i=v_i$ true in $T$.

For example, solvability of the equation $a=x_1bx_2$ ($a,b\in S$) means that $b$ divides $a$ in a finite semigroup containing $S$. The problem of eventual solvability of that equation was known as the \emph{Rhodes problem} since the 60s. Similar problems were known to be decidable  (Lyapin \cite{Ly}). Rhodes' problem was believed to be decidable also, and some partial results in that direction were proved. For example, Hall and Putcha \cite{HP} proved that the solvability of that equation over $S$ is decidable provided $a$ and $b$ are regular elements. Nevertheless, using Slobodskoi's result and split systems, we proved

\begin{theorem}[Kublanovsky, Sapir \cite{KuSRhodes}] (i) There is no algorithm to decide, given a finite 4-nilpotent
semigroup $S$ and two elements $a, b$ in $S$, whether there exists a
bigger finite semigroup $T > S$ such that $a = x_1bx_2$ for some $x_1,
x_2\in T$.

(ii) There is no algorithms to decide, given a finite 4-nilpotent semigroup
 $S$ and two elements $a, b$ in $S$, whether there exists a bigger finite
semigroup  $T > S$ such that $a = x_1bx_2$ and $b = x'_1ax_2'$ for some $x_1, x_1', x_2, x_2'\in T$.
\end{theorem}

Similar problems turned out to be undecidable for finite associative rings as well.

\subsubsection{Profinite groups and the restricted Burnside problem for general algebras} Recently Theorem \ref{t:sl} found several interesting and unexpected applications in the theory of profinite groups. In particular, Bridson and Wilton proved

\begin{theorem}[Bridson, Wilton \cite{BW}]  There are recursive sequences of finite presentations for residually finite groups $G_n$ and $\Gamma_n$ with explicit monomorphisms $u_n: G_n\hookrightarrow \Gamma_n$ such that (1) the profinite completions $\hat G_n$ and $\hat \Gamma_n$ are isomorphic  if and only if the induced map $\hat u_n$ is an isomorphism; (2) $\hat u_n$ is an isomorphism if and only if $u_n$ is surjective; and (3) the set $\{n\in \N\mid \hat G_n \hbox{ is not isomorphic to } \hat \Gamma_n\}$ is recursively enumerable but not recursive.
\end{theorem}

They also applied Slobodskoi's theorem to prove

\begin{theorem}[Bridson, Wilton \cite{Bridson}]\label{t:bw} There is no algorithm which, given a finite presentation of a group $G$ decides if $G$ has only finitely many finite quotients.
\end{theorem}

That result, in turn, was used in \cite{SapirBook} to prove that there is no algorithm which, given a finite set of identities of general algebras, decides whether the restricted Burnside property is true in the variety defined by these identities. More precisely we proved

\begin{theorem}[Sapir, \cite{SapirBook}]  There is no algorithm to decide, given a finite set of
identities of some type, whether for every $n$ there are only finitely many finite $n$-generated algebras satisfying these identities.
\end{theorem}

\section{Other algorithmic applications of Minsky machines}

Historically the first application of Minsky (register) machines was the result by Davis, Putnam and Robinson \cite{DPR} of unsolvability exponential diophantine equations: they associated an exponential diophantine equation $D(M)$ to every Minsky machine and its configuration $(1;m,0)$, and then proved that $M$ accepts the configuration $(1;m,0)$ if and only if $D(M)$ has an integer solution. This result turned out to be the first major step in solving Hilbert's 10th problem.

\subsection{Collatz type problems} A probably lesser known application of Minsky machines is the John Conway's version of Collatz problem \cite{Conway}. Recall that Collatz problem concerns with the function $\kappa\colon \N\to \N$ that takes every even number $n$ to $n/2$ and every odd number $n$ to $3n+1$. The problem is whether for every number $n$, $\kappa^s(n)=1$ for some $s$. Conway generalised this problem as follows. Let $\kappa$ be a piece-wise linear function $\N\to \N$ with finitely many pieces. Can we decide, given $n\in \N$, whether $\kappa^s(n)=1$ for some $s$?

Conway \cite{Conway} showed that the answer is negative even if we assume that the linear functions are just dilations (i.e., linear functions without the translation part). A simulation of Minsky machines is in that case not difficult. Let $M$ be a 2-glass Minsky machine with commands $\#\# 1,2,...,N,0$. We encode every configuration $(i;m,n)$ of $M$ by the number $p_i2^m3^n$ where $p_i$ is the $i+3$'d prime (that is $p_0=5, p_1=7$, etc.). Commands of $M$ are then encoded as pieces of the piece-wise dilation function $\kappa$:

$$
\begin{array}{|l|l|l|}
\hline
\hbox{Command \# i of } M& \hbox{value of }\kappa(n) & \hbox{condition on } n \\
\hline
\Add(1);j & \frac{2p_j}{p_i} n & p_i \hbox{ divides } n\\
\hline
\Add(2);j & \frac{3p_j}{p_i} n & p_i \hbox{ divides } n\\
\hline
\epsilon_{1}> 0 \to \Sub(1); j & \frac{p_j}{2p_i}n & 2p_i \hbox{ divides } n\\
\hline
\epsilon_{2}>0 \to \Sub(2);j & \frac{p_j}{3p_i}&2p_i \hbox{ divides } n\\
\hline
\epsilon_1=0\to j & \frac{p_j}{p_i} n & p_i \hbox{ divides } n \hbox{ but } 2\hbox{ does not divide } n\\
\hline
\epsilon_2=0\to j & \frac{p_j}{p_i} n & p_i \hbox{ divides } n \hbox{ but } 3\hbox{ does not divide } n\\
\hline
\hbox{Stop } (i=0) & \frac{1}{5} n& n=5\\
\hline
\end{array}
$$

In all other cases (say, when none of $p_i$ divides $n$), we set $\kappa(n)=n$.  Then it is easy to show that for $n=7\cdot 2^m$ there exists $s$ such that $\kappa^s(n)=1$ if and only if the configuration $(1;m,0)$ is accepted by $M$, and Conway's statement follows from Theorem \ref{t:MM}.

\subsection{Amalgams of finite semigroups}

Let $D$ and $E$ be two semigroups generated by sets $X$ and $Y$ respectively. Let $U=X\cap Y$, and the subsemigroups generated by $U$ in $D$ and in $E$ coincide. Then we say that $D\cup \la U\ra E$ is an amalgam of the semigroups $D$ and $E$ with amalamated subsemigroup $\la U\ra$. The corresponding amalgamated free product is the semigroup defined by the generating set $X\cup Y$ and all relations of $D$ and $E$.

It is well known that any amalgam of two finite groups is
embeddable into a group \cite{LS}. Moreover the free product with
amalgamation of two finite groups and in general the fundamental
group of any graph of finite groups has a free subgroup of finite index
(Karras, Pietrowski, Solitar, \cite{KPS}) so it is residually
finite. Hence it has solvable word problem and any amalgam of finite groups is embeddable into
a finite group. The situation with semigroup amalgams is quite different,  and there are many papers that tried to clarify the situation (see the introduction of \cite{SapirAmal}). In some sense these efforts were finalized by two results from \cite{SapirAmal} (very similar proofs of Parts (b), and (c)  were obtained independently and almost simultaneously by Jacson \cite{Jac}):

\begin{theorem}[Sapir, \cite{SapirAmal}]\label{th1} (a)
There exists an amalgam of two finite semigroups such that the
word problem is undecidable in the corresponding free product with
amalgamation and the amalgam embeds in the amalgamated free product.

(b) The problem of whether an amalgam of two finite
semigroups is embeddable into a semigroup is undecidable.

(c) The problem of whether an amalgam of two finite
semigroups is embeddable into a finite semigroup is undecidable.
\end{theorem}

The proofs of Parts (b), (c) uses Slobodskoi's Theorem \ref{t:sl} (and hence, implicitly, Minsky machines). The proof of Part (a) of Theorem \ref{th1} provides a simulation of a Minsky machine $M$ with $N+1$ commands $\#\# 1,....,N,0$ in an amalgamated product of two finite (even nilpotent) semigroups $D(M)$ and $E(M)$ described as follows.

First let us describe the generating sets of $D(M)$ and $E(M)$. These generating sets intersect by the set $U(M)$ which consists of $0, q_0, q_1, u_{i,j}$, $i=0,...,N$, $j=1,2$, $0$ acts as zero in both $D(M)$ and $E(M)$.
The semigroup $D(M)$ is generated by the
union of the set $\{a, \bb, q_i, p_i\ |\ i=0,\ldots, N\}$
and the
set $U$; $E(M)$ is generated by the
union of the set
$\{A, b, \ab, B\}$ and  $U$.

We shall see that the set $U(M)$ is in fact a subsemigroup with zero product: the product of every two elements in $U(M)$ is equal to $0$. The other relations of $D(M)$ and $E(M)$ are not as transparent and are much less transparent than in the cases considered above, so we start with explaining how these relations simulate the Minsky machine.

The word in the amalgamated free product $D(M)*_{U(M)} E(M)$ that corresponds to a configuration $(i;m,n)$ of the machine $M$ is
$W(i;m,n)= A(ab)^mq_i(\ab\bb)^nB$. Thus $A,B$ correspond to the bottoms of the two glasses, $(ab)^m$ simulates the  coins in the first glass, $(\ab\bb)^n$ simulates the coins in the second glass.

Suppose that the $i$th command of $M$ adds a coin in the first glass. Then first we replace $q_i$ by $au_{i,1}p_i$ using a relation of $D(M)$. Then using a relation of $E(M)$ we replace $u_{i,1}$ by $bu_{i,2}$. As a result the word $A(ab)^mq_i(\ab\bb)^nB$ transforms into $A(ab)^{m+1}u_{i,2}p_i(\ab\bb)^nB$. Finally using a relation of $D(M)$ we replace $u_{i,2}p_i$ by $q_j$, and we produce the word $W(j;m+1,n)$ as desired. Other commands of $M$ are simulated in a similar manner.

Here is the list of Minsky relations of semigroups $D(M), E(M)$:

$$
\begin{array}{|l|l|l|}
\hline
\hbox{Command \# i of } M& \hbox{ relation in } D(M) & \hbox{relation in } E(M)\\
\hline
\Add(1);j & q_i = au_{i,1}p_i,   & u_{i,1}=bu_{i,2} \\
          & u_{i,2}p_i=q_j &\\
\hline
\Add(2);j & q_i = p_iu_{i,1}\ab   & u_{i,1}=u_{i,2}\bb \\
          &p_iu_{i,2}=q_j &                            \\
\hline

\epsilon_{1}> 0 \to \Sub(1); j & q_i = u_{i,1}p_i   & bu_{i,1}=u_{i,2} \\
                               & au_{i,2}p_i=q_j    &                  \\

\hline
\epsilon_{2}>0 \to \Sub(2);j &q_i = p_i u_{i,1}  & u_{i,1}\ab=u_{i,2}\\
                             &p_iu_{i,2}\bb=q_j &\\
\hline
\epsilon_1=0\to j & q_i = u_{i,1}p_i & Au_{i,1}=Au_{1,2}\\
                  & u_{1,2}p_i=q_j &                    \\
\hline
\epsilon_2=0\to j & q_i=p_iu_{i,1} & u_{i,1}B=u_{i,2}B\\
                  & p_iu_{i,2}=q_j &\\
\hline
\end{array}
$$

We also add to $D(M)$ (resp. $E(M)$) all relations of the form $w=0$ where $w$
is any word
in generators of $D(M)$ (resp. $T(M)$) which is not a subword of any word
participating in the Minsky relations, or the words $Aq_1B$ and $Aq_0B$.
For example, $a\ab=aq_i=q_i\ab=0$ in
$D(M)$, $Ab=AB=0$ in $E(M)$. Thus, in particular, $U(M)$ is indeed
the semigroup with zero product, and each semigroup $D(M)$ and $E(M)$ are 4-nilpotent and finite.

The proof that the amalgamated product $R$ of $D(M)$ and $E(M)$ simulates $M$ proceeds as follow. First we prove that if two configurations $(i;m,n)$ and $(i';m',n')$ of $M$ are equivalent, then the corresponding words $W(i;m,n)$ and $W(i';m',n')$ are equal in $R$. Then we notice that the presentation of $R$ is confluent \cite{SapirBook}. Hence if two words $W, W'$ are equal in $R$, then there exist two  sequences of applications of relations of $R$ from left to right, one starting at $W$, another starting at $W'$ which end at the same word $W''$. This implies that if $W(i;m,n)=W(i';m',n')$ in $R$, then the configurations $(i;m,n)$ and $(i';m',n')$ of the Minsky machine $M$ are equivalent.

\subsection{Complicated residually finite semigroups}

It is well known that finitely presented residually finite algebras (of finite signature) are much simpler algorithmically than arbitrary finitely presented algebras. For example, the word problem in every such algebra is decidable (see McKinsey's algorithm in \cite{Malcev}). Moreover the most ``common" residually finite algebras, say, the linear groups over fields, are algorithmically ``tame": the word problem in any linear group is decidable in polynomial time and even log-space \cite{LZ77}. Surprisingly till \cite{KMS} not much was known about possible complexity of arbitrary finitely presented residually finite algebras, even in the cases of semigroups and groups.

Recall that the McKinsey's algorithm for solving the word problem in a finitely presented algebra $A$ consists of two competing parts running in parallel. The first part enumerates all finite homomorphic factors of $A$ and checks if the images of two elements $a,b$ are different in one of these factors. The second part enumerates all consequences of the defining relations of $A$. If $a\ne b$ in $A$, then the first part will eventually ``win", if $a=b$, then the second part ``wins". In either case we will eventually know whether or not $a=b$ in $A$.

Thus there are three ways to estimate the complexity of a residually finite algebra:

\begin{enumerate}

\item[(1)] The computational complexity of the word problem.

\item[(2)] The \emph{depth function} $\rho_A(n)$ which is the smallest function $\N\to \N$ such that  given two words $u, v$ of size $\le n$ \footnote{By a ``word" I understand any term involving operations of $A$. It is an ordinary word in the case of semigroups, a word possibly containing inverses of generators in the case of groups, a non-commutative polynomial in the case of rings, etc. The size of a word - the number of symbols needed to write it.} in generators of $A$, such that $u\ne v$ in $A$, there exists a homomorphism from $A$ onto a finite algebra of size at most $\rho_A(n)$ which separates $u$ and $v$.

\item[(3)] The \emph{Dehn function} $d_A(n)$ which is the smallest function $\N\to \N$ such that, given two words $u, v$ of size $\le n$ in generators of $A$, such that $u=v$ in $A$ one needs at most $d_A(n)$ applications of relations of $A$ to deduce the equality $u=v$.
\end{enumerate}

Note that the Dehn function and the depth function are two of the most important asymptotic characterics of an algebra. Both functions are recursive for every finitely presented residually finite algebra (of finite signature). Gersten \cite{Ger, GR} asked for a bound of the Dehn function of a finitely presented linear group. The answer is still not known. The depth function of groups was first studied  by Bou-Rabee \cite{Bou-Rabee}. It is known (and easy) that for linear groups, it is at most polynomial. In fact till \cite{KMS}, no finitely presented group or semigroup with Dehn function or depth function greater than an exponent was known.

The next result from \cite{KMS} shows that there are finitely presented residually finite semigroups and groups with arbitrary high complexity in each of the three ways to measure the complexity, and also semigroups and groups with word problem in $P$ but arbitrary large (recursive) Dehn and depth functions.

\begin{theorem} \label{t:kms} Let $f\colon \N\to \N$ be any recursive function. Then

(i) there is a  residually finite finitely presented group that is
solvable of class 3  with Dehn function $d_G$ such that $d_G\succcurlyeq f$. In addition,  one can make the group $G$ such that the time complexity of the word problem in  $G$ is at least as large  as any given recursive function or one can make $G$ such that the word problem is  in polynomial time.

(ii) there is a residually finite finitely presented
solvable of class 3 group $G$ with depth function greater than $f$. In addition, one can make the group $G$ such that the word problem in  $G$ is   at least as hard  as the membership problem in a given recursive set of natural numbers $\Z$ or one can make $G$ such that the word problem is in  polynomial time.
\end{theorem}

Here for two functions $f, g\colon \N\to \N$ we write $f\succcurlyeq g$ if for some constants $c_1, c_2,c_3\ge 0$, we have $f(n)\ge c_1g(c_2n)-c_3$ for all $n$.

To illustrate the proof of Theorem \ref{t:kms}, we will present a construction of complicated residually finite semigroups. In the case of groups, the construction is based on a similar idea but is technically much more complicated.

We start with a deterministic Turing machine $TM$ that recognizes a recursive set $Z$ of natural numbers. By \cite{MDavis}, we can assume that $TM$ is universally halting.  This means that $TM$ has only finite number of computations starting from any given configuration. We modify $TM$ to obtain a $\Sym$-universally halting machine $TM'$. That means that the (non-deterministic) Turing machine $\Sym(TM')$ is universally halting. The construction from \cite{malcev} (used to prove Theorem \ref{t:MM} above) produces then a $\Sym$-universally halting 3-glass Minsky machine $M_3$ that recognizes the set $Z$. One can easily modify the construction of the semigroup $\orr S_2$ above to produce a finitely presented semigroup $\orr S_3$ that simulates the 3-glass Minsky machine $M_3$. As above every word in the generators of $\orr S_3$ which is non-zero in $\orr S_3$ is a subword of a word $W_3(i;m,n,k)$ corresponding to a configuration $(i; m,n,k)$ of $M_3$ modulo the glass permuting relations. The semigroup $\orr S_3$ is residually finite. Indeed, take any two words $u, v$ in the generators of $\orr S_3$, such that $u\ne v$ in $\orr S_3$. We can assume that $u\ne 0$ in $\orr S_3$. Then since $M_3$ is $\Sym$-universally halting, there are only finitely many elements of $\orr S_3$ that divide $u$. It is not difficult to prove that if $v$ divides $u$, then $v\ne 0$ in $\orr S_3$ and $u$ does not divide $v$. Thus in that case we can replace $v$ for $u$. So we can assume that $v$ does not divide $u$ in $\orr S_3$. As above take the ideal $J$ of all elements in $\orr S_3$ that do not divide $u$. Then the Rees quotient $\orr S_3/J$ is finite and $u\ne v$ in the quotient which proves residual finiteness. 

Now the word problem in $\orr S_3$ polynomially reduces to the configuration equivalence problem for $M_3$. It is at least as hard as the membership problem for $Z$, and in fact by carefully choosing $TM'$, we can make that problem polynomially equivalent to the membership problem for $Z$. Therefore we can make the complexity of the word problem in $\orr S_3$ as hard or as easy (i.e., at most polynomial time) as we want. The Dehn function of $\orr S_3$ is at least as large as the time function of $TM'$ (in fact much larger). Even in the case when $Z$ is in $P$, we can construct a Turing machine recognizing $Z$ and with arbitrary large (recursive) time function. Say, after the machine wants to stop, we make it compute something really complicated, and only then stop. This way we construct a finitely presented residually finite semigroup with polynomial time complexity of the word problem and arbitrary high recursive Dehn function.

To make the depth function high, we add two more glasses to the Minsky machine $M_3$ and modify the commands to obtain a 5-glass non-deterministic Minsky machine $M_5$. We modify the commands of $M_3$ as follows. First to every command of $M_3$, we add the instruction to add a coin to glass $\# 4$ provided glass $\# 5$ is empty. Also assuming $M_3$ had $N+1$ command, for every $i=1,...,N$ we add two commands $\# i$:

\begin{equation}
\label{e:01}
i; (\Add(4), \Add(5)), i.
\end{equation}

and

\begin{equation}
\label{e:02}
i; (\epsilon_{4}=0, \epsilon_{5}=0)\to  0
\end{equation}

That is executing command $i$ the machine can add as many coins (equal amounts) to glasses $\# 4, \# 5$ and if both these glasses are empty, then the machine can stop. It is not difficult to prove, as above, that the semigroup $\orr S_5$ simulating $M_5$ is residually finite and its word problem is polynomially equivalent to the membership problem for $Z$. Let $a_4$ and $a_5$ be the generators of $\orr S_5$ simulating coins in glasses $4$ and $5$. In order to ensure that the depth function of $\orr S_5$ is high we take a word $W_5=W_5(1;z,0,0,0,0)$ of length $n$, corresponding to the input configuration such that $z\not\in Z$ and there exists a very long computation, say of length $L\gg n$ of $M_5$ starting at $W_5$ and not using commands (\ref{e:01}) (we can always assume modify the Turing machine $TM'$ to make this happen). Note that since $z\not\in Z$, $W_5\ne 0$ in $\orr S_5$. Suppose that there exists a homomorphism from $\orr S_5$ to a finite semigroup $E$ of order $l$, and, say, with $l! < L$,  which separates $W_5$ from $0$. Then for every element $x$ in $E$ we will have $x^m=x^{2m}$ for $m=l!$. Consider a very long computation of $M_5$ starting with the configuration $(1;z,0,0,0,0)$. In that computation, we must get a configuration $\delta$ where the glass $\# 4$ contains $m$ coins and glass $\# 5$ empty. The corresponding word $W_5(\delta)$ will contain the subword $a_4^mA_4A_5$. Applying now relations corresponding to the command (\ref{e:01}), we can change this subword (without touching the rest of the word) to $a_4^{2m}A_4a_5^mA_5=a_4^mA_4a_5^mA_5$ in $E$. Applying the relations corresponding to $(\ref{e:01})$ again (this time from right to left), we obtain a word that is equal to $W_5$ in $E$ and has subword $A_3A_4A_5$ (it corresponds to a configuration with empty glasses $\#\# 4,5$). Applying now the relations corresponding to (\ref{e:02}), and then the relation $q_0=0$, we deduce that $W_5=0$ in $E$, a contradiction. This contradiction shows that the factorial of the order $|E|$ cannot be smaller than $L$. Thus the depth function of $\orr S_5$ can be as large as we want.

\end{document}